\documentclass{article}
\date{}
\usepackage{hyperref}
\usepackage{indentfirst}
\usepackage{amsmath}
\usepackage{amssymb}
\usepackage{amsthm}
\title{Two dimensional disjoint minimal graphs}
\author{Linfeng Zhou}
\begin{document}
\maketitle

\begin{abstract} In this paper, under the assumption of Gauss curvature vanishing at infinity, we will prove Meeks' conjecture: the number of disjointly supported minimal graphs in $\mathbb{R}^3$ is at most two.   
\\

\end{abstract}

\newtheorem{Th}{Theorem}[section]
\newtheorem{Prop}[Th]{Propositon}
\newtheorem{Col}[Th]{Corollary}
\newtheorem{Lem}[Th]{Lemma}
\newtheorem{Ex}[Th]{Example}
\newtheorem{Con}[Th]{Conjecture}
\newtheorem*{Claim}{Claim}
\newtheorem*{Remark}{Remark}

\section{Introduction}
Let $\Omega$ be an open subset in $\mathbb{R}^2$ and denote its boundary by $\partial\Omega$. As we know, if a function $u(x)$ which is defined on $\Omega$ satisfies the equation 
\begin{equation}\label{ms}div(\frac{\nabla u}{\sqrt{1+|\nabla u|^2}})=0,\end{equation}
$G=\{(x,u(x))|x\in \Omega\}$ is called a minimal graph in $\mathbb{R}^3$. Furthermore, we call the minimal graph $G$ is supported on $\Omega$ if $u|_{\partial \Omega}=0$ and $u\geq 0$.

Meeks \cite{M} has conjectured that the number of disjointly supported minimal graphs with zero boundary values over an open subset in $\mathbb{R}^2$ is at most 2. In fact, for arbitrary dimension, Meeks-Rosenberg \cite{M-R} proved if a set  of disjointly supported minimal graphs have bounded gradient, then the number of the graphs must be finite. Later, Li-Wang \cite{L-W2} gave an upper bound of the number of the graphs without any assumption on the growth rate of each graph. As a corollary, when minimal graphs are two dimensional in $\mathbb{R}^3$, they obtained the number is at most 24.   At the same time, Spruck \cite{S} proved that there are at most two admissible sub-linear growth solution pairs of the equation (\ref{ms}) defined over disjoint domains. Recently, by using angular density, Tkachev \cite{T} showed the number of two dimensional disjointly supported minimal graphs  is less than or equals 3.  

Observing the similarity between the disjoint d-massive set and disjointly supported minimal graphs,  we can apply the method of proving finiteness theorem of disjoint $d$-massive sets in $\mathbb{R}^2$ \cite{L-W1} to study disjoint minimal graphs.  Actually, we obtain the following theorem:   

\begin{Th} \label{mt}Suppose $\{G_i=(\Omega_i,u_i)\}_{i=1}^k$ is a set of disjointly supported minimal graphs in $\mathbb{R}^3$ where each $\Omega_i$ is an open subset  in $\mathbb{R}^2$. If the Gauss curvature $K_i(x)$ of each graph satisfies \[K_i(x)\rightarrow 0\qquad (|x|\rightarrow\infty),\] then the number of the graphs $k$ is at most two.
\end{Th}

By choosing a slight different region of integration, one has a stronger result comparing  with the theorem of Spruck\cite{S}.
\begin{Col}\label{st}
Suppose $\{G_i=(\Omega_i,u_i)\}_{i=1}^k$ is a set of disjointly supported minimal graphs in $\mathbb{R}^3$ where each $\Omega_i$ is an open subset  in $\mathbb{R}^2$. If each graph has sub-linear growth, then $k$ is at most two.
\end{Col}

 This work was carried out during the author visiting University of Minnesota. He would like to express his gratitude to Professor Jiaping Wang for his guidance and many helpful discussion. The author also appreciate the mathematical department of UMN for its hospitality. Thanks for the referees for many valuable comments. 
\section{Proof of theorem \ref{mt} }

In the following,  we denote the 3-dimensional ball of radius $R$ centered at the origin of $\mathbb{R}^3$ by $B^3(R)$  and the 2-dimensional sphere of radius $R$ by $S^2(R)$. Actually, the key is to establish a refined estimate of the sum of all curves' length $\ell( G_i\cap S^2(R))$ 
when $R$ is sufficiently large.   

\begin{Th}\label{le}Suppose $\{G_i=(\Omega_i,u_i)\}_{i=1}^k$ is a set of disjointly supported minimal graphs in $\mathbb{R}^3$ where the Gauss curvature $K_i(x)$ of each $G_i$ satisfies \[K_i(x)\rightarrow 0\qquad (|x|\rightarrow\infty).\] For a sufficiently large radius $R$, $\sum\limits_{i=1}^k\ell(G_i\cap S^2(R))$ is bounded by
\[\sum_{i=1}^k\ell(G_i\cap S^2(R))\leq \pi^2R+o(1)R.\]
In particular, when $k=3$, we have an refined estimate
\[\sum_{i=1}^3\ell(G_i\cap S^2(R))\leq 2\sqrt{2}\pi R+o(1)R.\]
\end{Th}

Before proving the theorem \ref{le}, we introduce a lemma. 
\begin{Lem} \label{ap} Let $B^3_+(R)$ be a 3-dimensional upper half ball with the radius $R$  and $S^2_+(R)$ be a 2-dimensional upper half sphere. Suppose $\pi_i:G_i\rightarrow \mathbb{R}^2$ is  the natural projective map.
If $\Sigma_1$, $\Sigma_2$, $\dots$, $\Sigma_s$ are the planes in $\mathbb{R}^3$ such that each interior of $\pi_i(\Sigma_i\cap B^3_+(R))$ does not intersect for a sufficiently large $R$, then the length of the curve $\Sigma_i\cap S^2_+(R)$ satisfies 
\[\sum_{i=1}^s\ell(\Sigma_i\cap S^2_+(R))\leq \pi^2 R.\]
Moreover, when $s=3$, we have a refined estimate
\[\sum_{i=1}^3\ell(\Sigma_i\cap S^2_+(R))\leq 2\sqrt{2}\pi R.\]
\end{Lem}

\begin{proof} Suppose $D(R)=\{(x_1,x_1,0)| x_1^2+x_2^2\leq R^2\}$ is a disk in $\mathbb{R}^3$. Since each $\Sigma_i$  is a plane, $\Sigma_i\cap D(R)$ is a chord and let $\theta_i$ be the corresponding central angle. 
Here we only need to consider the case that the union of each chords $\cup_{i=1}^s(\Sigma_i\cap D(R))$ is a polygon.  Otherwise, one can add more planes which still satisfy the required conditions such that above intersection yields a polygon.  

If the centre of the disk $D(R)$ is in the interior of the polygon or on one of the edge of the polygon, this means each central angle $\theta_i$ satisfies $0<\theta_i\leq \pi$.
For each interior of  $\pi_i(\Sigma_i\cap B^3_+(R))$ does not intersect, by a simple computation one can obtain the following inequality about the length of the arc $\ell(\Sigma_i\cap S^2_+(R))$ 
\[\ell(\Sigma_i\cap S^2_+(R))\leq\pi R \sin\frac{\theta_i}{2}.\] 
The RHS achieves the maximum if and only if $\Sigma_i$ is perpendicular to the disk $D(R)$. Thus
\begin{eqnarray}\label{leiq1}
&\sum\limits_{i=1}^s\ell(\Sigma_i\cap S^2_+(R))&\leq \sum_{i=1}^s \pi R\sin\frac{\theta_i}{2}\leq \pi R s \sin(\frac{1}{s}\sum_{i}^s \frac{\theta_i}{2})\nonumber\\
&&\leq \pi R s \sin(\frac{\pi}{s})\leq \pi^2 R.
\end{eqnarray} In the second above inequality, we use the concave property of the sine function on the interval $[0,\pi]$. 

For a special case when $s=3$, from (\ref{leiq1}) one can yield
\begin{equation}\label{leiqc1}\sum\limits_{i=1}^3\ell(\Sigma_i\cap S^2_+(R))\leq 3\pi R \sin(\frac{\pi}{3})=\frac{3\sqrt{3}}{2}\pi R.\end{equation}

If the centre of the disk $D(R)$ is outside the polygon, namely there exists an $i_0$ such that $\theta_{i_0}>\pi$. For simplicity, let us assume $i_0=s$. A similar computation induces that 
\[\ell(\Sigma_i\cap S^2_+(R))\leq\pi R \sin\frac{\theta_i}{2}\quad \text{for}\quad 1\leq i\leq s-1,\]
\[\ell(\Sigma_s\cap S^2_+(R))\leq R\theta_s.\]
The first equality holds if and only if $\Sigma_i$ is perpendicular to the disk and the second equality holds if and only if $\Sigma_s$ is in the same plane of the disk $D(R)$.  Hence one will have
\begin{eqnarray}\label{leiq2}
& \sum\limits_{i=1}^s\ell(\Sigma_i\cap S^2_+(R))&\leq \sum_{i=1}^{s-1} \pi R\sin\frac{\theta_i}{2}+ R\theta_s\leq \sum_{i=1}^{s-1}\pi R\sin\frac{\theta_i}{2}+2\pi R \sin\frac{\theta_s}{4} \nonumber\\
&&\leq \pi R (s+1) \sin \frac{\pi}{s+1}\leq \pi^2 R.
\end{eqnarray}

If $s=3$, by (\ref{leiq2}) we obtain that
\begin{equation}\label{leiqc2}
 \sum\limits_{i=1}^3\ell(\Sigma_i\cap S^2_+(R))\leq  4\pi R \sin(\frac{\pi}{4})=2\sqrt{2}\pi R
\end{equation}

The conclusion is derived from (\ref{leiq1}), (\ref{leiq2}) and (\ref{leiqc1}), (\ref{leiqc2}).
\end{proof}

\begin{proof}[Proof of the theorem \ref{le}] For each minimal graph $G_i$,  since the Gauss curvature $K_i=0$ at infinity, it means $G_i$ is asymptotic to a flat plane. Therefore, we can use the intersection of a plane $\Sigma_i$ and $S_+^2(R)$ to approximate  the curve $G_i\cap S^2(R)$. By the lemma \ref{ap}, one has
\[\ell(G_i\cap S^2(R))\leq \ell(\Sigma_i\cap S_+^2(R))+o(1)R.\]
Therefore 
\[\sum_{i=1}^k\ell(G_i\cap S^2(R))\leq \sum_{i=1}^k\ell(\Sigma_i\cap S_+^2(R))+o(1)R\leq \pi^2R+0(1)R.\]
\end{proof}

The following lemma of the area growth estimate of a minimal graph is well-known argument, and one can see \cite{L-W2} for the details 

\begin{Lem} \label{Le1} Let  $G=(\Omega,u)$ be a minimal graph in $\mathbb{R}^3$, the area of $G\cap B^3(R)$ satisfies
\[A(G\cap B^3(R))\leq 3\pi R^2.\]
\end{Lem}

We are now ready to prove the main theorem.
\begin{proof}[Proof of the theorem \ref{mt}]Let $B^3(R)$ be the ball of radius $R$ in $\mathbb{R}^3$. Since
\[\int_{G_i\cap B^3(R)}|\tilde{\nabla} u_i|^2\leq \int_{G_i\cap \partial B^3(R))}u_i(\tilde{\nabla} u_i\cdot \frac{\partial }{\partial r})\]
where $\tilde{\nabla}$ means the gradient operator on $G_i$,
one has
\[2\lambda_1^{\frac{1}{2}}(G_i\cap\partial B^3(R))\int_{G_i\cap B^3(R))}|\tilde{\nabla}u_i|^2\leq 2\lambda_1^{\frac{1}{2}}\int_{G_i\cap \partial B^3(R)}u_i\cdot\frac{\partial u_i}{\partial r}\]
\[\leq \lambda_1\int_{G_i\cap \partial B^3(R)}u_i^2+\int_{G_i\cap\partial B^3(R)}(\frac{\partial u_i}{\partial r})^2\]
\[\leq\int_{G_i\cap \partial B^3(R)}|\bar{\nabla} u_i|^2+\int_{G_i\cap \partial B^3(R)}(\frac{\partial u_i}{\partial r})^2\]
\[=\int_{G_i\cap \partial B^3(R)}|\tilde{\nabla}u_i|^2.\]
Here $\lambda_1^{\frac{1}{2}}(G_i\cap\partial B^3(R))$ denotes the first Dirichlet eigenvalue
on $G_i\cap\partial B^3(R)$. As we know, in $\mathbb{R}^3$ the following inequality holds:
\[\lambda_1^{\frac{1}{2}}(G_i\cap\partial B^3(R))\geq \frac{\pi^2}{\ell^2(G_i\cap\partial B^3(R))}.\]
Therefore
\[\frac{\int_{G_i\cap\partial B^3(R)}|\tilde{\nabla} u_i|^2}{\int_{G_i\cap B^3(R)}|\tilde{\nabla} u_i|^2}\geq 2\lambda_1^{\frac{1}{2}}\geq\frac{2\pi}{\ell(\Gamma_i)},\]
where $\Gamma_i:=G_i\cap\partial B^3(R).$
Thus we obtain
\[\sum_{i=1}^{k}\frac{\int_{G_i\cap\partial B^3(R)}|\tilde{\nabla} u_i|^2}{\int_{G_i\cap B^3(R)}|\tilde{\nabla} u_i|^2}\geq \sum_{i=1}^{k}\frac{2\pi}{\ell(\Gamma_i)}.\]
Notice that
\[k^2\leq (\sum_{i=1}^k \ell(\Gamma_i))(\sum_{i=1}^k\frac{1}{\ell(\Gamma_i)}).\]
According to the theorem \ref{le}, one has
\[\sum_{i=1}^k\ell(\Gamma_i)\leq \pi^2 R+o(1)R\]
for a sufficiently large radius $R$. 
Then it can be concluded
\begin{equation}\label{peq1}\sum_{i=1}^{k}\frac{\int_{G_i\cap\partial B^3(R)}|\tilde{\nabla} u_i|^2}{\int_{G_i\cap B^3(R)}|\tilde{\nabla} u_i|^2}\geq \frac{2\pi k^2}{R(\pi^2+o(1))}.\end{equation}
Observing that
\begin{equation}\label{peq2}\int_{G_i\cap \partial B^3(r)}|\tilde{\nabla} u_i|^2=\frac{\partial}{\partial r}\int_{G_i\cap B^3(r)}|\tilde{\nabla}u_i|^2.\end{equation}
From (\ref{peq1}) and (\ref{peq2}) to obtain
\begin{equation}\label{peq3}\ln(\prod_{i=1}^k \frac{\int_{G_i\cap B^3(R)}|\tilde{\nabla} u_i|^2}{\int_{G_i\cap B^3(R_0)}|\tilde{\nabla} u_i|^2} )\geq \frac{2\pi k^2}{\pi^2+o(1)}\ln(\frac{R}{R_0}).\end{equation}
On the other hand, let $(x,y,u_i(x,y))$ be a parametrization of $G_i$, then the induced metric on $G_i$ is 
\[ds^2=(1+(u_i)_x^2)dx^2+2(u_i)_x(u_i)_ydxdy+(1+(u_i)_y^2)dy^2.\]
Hence
\[|\tilde{\nabla} u_i|=\sqrt{u_{x^i}u_{x^j}g^{ij}}=\sqrt{\frac{|\nabla u_i|^2}{1+|\nabla u_i|^2}}\leq 1.\]
From it one can deduce
\begin{equation}\prod_{i=1}^k\int_{G_i\cap B^3(R))}|\tilde{\nabla} u_i|^2\leq A^k(G_i\cap B^3(R))\label{peq4}\leq (3\pi R^2)^k.\end{equation}
Combining (\ref{peq3}) and (\ref{peq4}) implies
\[\frac{2\pi k^2}{\pi^2+o(1)}(\ln R-\ln R_0) \leq 2k\ln R+c_1 .\]
Let $R\rightarrow +\infty$ to have
\[k\leq \pi.\]
This inequality indicates that $k\leq 3$.

If $k=3$, repeating the above process and using the refined length estimate in the theorem \ref{le} provides
\[k\leq 2\sqrt{2}\]
which is a contradiction. 

Thus $k$ has to be at most 2. 
 \end{proof}

\begin{Remark} In \cite{T}, Tkachev has already proved the number of two dimensional disjointly supported minimal graphs is at most 3. Here a different approach can lead to a better estimate if assuming the Gauss curvature vanishes at infinity. 
\end{Remark}

\section{Proof of corollary \ref{st} }
Let $\pi_i:G_i\rightarrow \mathbb{R}^2$ be the natural projective map and $B^2(R)$ be the ball of radius $R$ in $\mathbb{R}^2$. By employing the same method in the proof of theorem \ref{mt} except for using a different region of integration $\pi_i^{-1}(\Omega_i\cap B^2(R))$, one can conclude

 \begin{Th} Suppose $\{(\Omega_i,u_i)\}_{i=1}^k$ is a set of disjointly supported minimal graphs in $R^3$ where each $\Omega_i$ is an open subset in $R^2$. If the gradient of each $u_i$ is bounded by $c$, i.e. $|\nabla u_i|\leq c$ , then the number $k\leq 2\sqrt{1+c^2}$. 
\end{Th}

\begin{proof}: 
By a similar argument,  one can obtain that
\[\sum_{i=1}^{k}\frac{\int_{\pi_i^{-1}(\Omega_i\cap\partial B^2(R))}|\tilde{\nabla} u_i|^2}{\int_{\pi_i^{-1}(\Omega_i\cap B^2(R))}|\tilde{\nabla} u_i|^2}\geq \frac{2\pi k^2}{\sum\limits_{i=1}^{k}\ell(\Gamma_i)}.\]
where $\Gamma_i:=\pi_i^{-1}(\Omega_i\cap\partial B^2(R)).$
If one chooses the parameter $(R\cos(\theta),R\sin(\theta),u_i(R\cos(\theta),R\sin(\theta)))$ of the curve $\Gamma_i$ and assume $|\nabla u_i|\leq c$, then
\[\ell(\Gamma_i)=\int_{\theta_0}^{\theta_1}\sqrt{R^2+[-(u_i)_xR\sin(\theta)+(u_i)_yR\cos(\theta)]^2}d\theta\]
\[\leq \int_{\theta_0}^{\theta_1}\sqrt{R^2+[(u_i)_x^2+(u_i)_y^2](R^2\sin(\theta)^2+R^2\cos(\theta)^2)}d\theta\]
\[\leq (\theta_1-\theta_0)R\sqrt{1+c^2}.\]
Since the minimal graphs are disjoint, so
\[\sum_{i=1}^k\ell(\Gamma_i)\leq 2\pi R\sqrt{1+c^2}.\]
Then it can be concluded
\begin{equation}\label{eq1}\sum_{i=1}^{k}\frac{\int_{\pi_i^{-1}(\Omega_i\cap\partial B^2(R))}|\tilde{\nabla} u_i|^2}{\int_{\pi_i^{-1}(\Omega_i\cap B^2(R))}|\tilde{\nabla} u_i|^2}\geq \frac{k^2}{R\sqrt{1+c^2}}.\end{equation}
 Integrating (\ref{eq1}), one will obtain
\begin{equation}\label{eq3}\ln(\prod_{i=1}^k \frac{\int_{\pi_i^{-1}(\Omega_i\cap\partial B^2(R))}|\tilde{\nabla} u_i|^2}{\int_{\pi_i^{-1}(\Omega_i\cap B^2(R_0))}|\tilde{\nabla} u_i|^2} )\geq \frac{k^2}{\sqrt{1+c^2}}\ln(\frac{R}{R_0}).\end{equation}
On the other hand, 
\begin{eqnarray}\label{eq4}
\prod_{i=1}^k\int_{\pi_i^{-1}(\Omega_i\cap B^2(R))}|\tilde{\nabla} u_i|^2&\leq& A^k(\pi_i^{-1}(\Omega_i\cap B^2(R)))=(\int_{\Omega_i\cap B^2(R)}\sqrt{1+|\nabla u|^2})^k\nonumber\\
&\leq& (\sqrt{1+c^2}\pi R^2)^k.
\end{eqnarray}
Combining (\ref{eq3}) and (\ref{eq4}), we have
\[\frac{k^2}{\sqrt{1+c^2}}(\ln R-\ln R_0) \leq 2k\ln R+c_1 .\]
Let $R\rightarrow +\infty$ to have
\[k\leq 2\sqrt{1+c^2}.\]
\end{proof}

Obviously, corollary \ref{st}  follows from above theorem when each graph satisfies 
\[|\nabla u_i|\rightarrow 0\quad (|x|\rightarrow +\infty).\]

\begin{Remark} J. Spruck has already proved the corollary \ref{st} under the assumption of a certain decay rate of  Gauss curvature at infinity \cite{S}.  However, here we do not need any kind of restrictions on Gauss curvature.   
\end{Remark}

\noindent\emph{Department of Mathematics\\East China Normal University\\Shanghai, 200241, China}\\

\noindent\emph{Department of Mathematics\\University of Minnesota\\ Minneapolis, MN, 55455, U.S.A.}\\

\noindent\emph{E-mail: lfzhou@math.ecnu.edu.cn}

 \LaTeX
\end{document}